\documentclass[letterpaper,12pt]{article}
\usepackage[margin=1in]{geometry}
\usepackage[T1]{fontenc}
\usepackage{lmodern,microtype}
\usepackage{xcolor}
\usepackage{enotez}
\setenotez{backref=true,list-name=Notes,list-heading={}}
\usepackage[colorlinks=true,urlcolor=blue]{hyperref}
\title{Will AI kill us all---or just make mathematicians obsolete?}
\author{Sergiu Klainerman}
\date{}
\let\footnote\endnote

\begin{document}
\maketitle

\begin{abstract}
AI's recent advances have fueled predictions of human intellectual obsolescence and fears for the future of mathematics. I argue that, impressive as these achievements are, they do not establish AI's capacity to invent new concepts and theories or sustain autonomous recursive self-improvement. Human intelligence also involves judgment, imagination, and embodied experience that standard measures fail to capture. Mathematics extends beyond proving theorems: its creative core lies in conceptual discovery and understanding. AI can strengthen this enterprise, but realizing its promise requires protecting attribution, open exchange, confidentiality, and mathematical training. I propose an academic AI consortium to provide powerful research tools with enforceable guarantees safeguarding unpublished ideas.
\end{abstract}

\subsection{AI will kill us all}

In 2023, hundreds of top AI researchers signed a
\href{https://www.safe.ai/statement-on-ai-risk}{one-sentence public
statement} putting the risk of extinction from AI on the same footing as
pandemics and nuclear war. 

\textit{Mitigating the risk of extinction from AI should be a global priority alongside other societal-scale risks such as pandemics and nuclear war}

Paradoxically, the  statement was signed by some of the most aggressive pursuers   of the technology, such as  the CEOs of  OpenAI and Anthropic,  arguably  the top AI companies in the world.

That warning became a rallying cry for
people working in AI who  envision  the development of  superintelligent   AI  (SAI) both as a Holy Grail to be pursued  and as a cool doomsday
scenario to scare all of us with. It was recently brought  to the center of national attention by Jacob Coxon, a researcher at Anthropic and
a former researcher at OpenAI, who accused  his former employers of racing towards self-improving  superintelligence, risking human extinction.  Coxon's own thread on X received more than 170 million views within days; 
 researchers at Anthropic backed his warning\footnote{  Anthropic's head of alignment, Evan Hubinger, publicly affirmed Coxon's
concern and stated that he personally thought there was a greater than
10\% chance of extinction \href{https://www.foxbusiness.com/technology/anthropic-researcher-says-ai-has-over-10-chance-kill-all-humans-within-next-decade}{within the next decade}. Other researchers offered higher personal estimates, under different assumptions: Geoffrey Irving, formerly at Google DeepMind and OpenAI, put the risk of human extinction at approximately 50\%, without specifying a time horizon in the cited statement (\href{https://frominside.ai/}{interview, 2:56}); Marcus Williams of OpenAI reportedly put it at 70\% within three years, conditional on the absence of regulation or a coordinated slowdown (\href{https://beincrypto.com/openai-researcher-human-extinction-odds/}{BeInCrypto, September 11, 2026}).}; the company's own CEO called for slowing the development of frontier
models; and panicked politicians   began talking about aggressive  federal oversight measures. It
does not take a prophet to see that the story has the potential to spark a
movement of its own, intensely vocal and moralizing,   organized on the model of climate
activism,  demanding not just  regulation but an outright halt to AI research.

In spite of the   internal warnings  of many  of its own researchers, including  its CEO, OpenAI kept racing  towards  SAI. On September  8, 2026,  the same day  when Coxon made public  his warnings, OpenAI  announced  an AI-generated solution to the Navier-Stokes problem, one of the seven Millennium Prize Problems.  The feat, said to have taken as many as 10,000 AI ``agents'' working in tandem and many millions of dollars' worth of computing, caps a dizzying run of ever more impressive achievements by artificial intelligence in mathematics (\url{https://arxiv.org/abs/2608.16753})--- the  discipline  that   represents  the  highest expression of  human rationality.   On October 6, OpenAI followed this announcement with a mass release of hundreds of AI-generated manuscripts claiming advances across a wide range of mathematical fields.\footnote{OpenAI, \href{https://openai.com/index/sharing-ai-progress-in-mathematics/}{``Sharing AI progress in mathematics,''} 6 October 2026; see also the accompanying \href{https://github.com/openai/math}{manuscript collection}.} OpenAI has already signaled further releases, deepening fears among mathematicians that their profession faces an existential threat.\footnote{The same OpenAI announcement explicitly refers to ``future releases.'' For accounts of the profession's fears, see Michel Schellekens, \href{https://proofsandprompts.com/2026/10/11/on-our-existential-crisis-and-trustworthy-ai/}{``On our existential crisis and trustworthy AI,''} \textit{Proofs and Prompts}, 11 October 2026; and David Larousserie, \href{https://www.lemonde.fr/en/science/article/2026/10/12/openai-s-latest-blow-to-mathematicians-they-have-annihilated-my-field-of-research_6758515_10.html}{``OpenAI's latest blow to mathematicians: `They have annihilated my field of research',''} \textit{Le Monde}, 12 October 2026.}

 If AI can solve mathematical problems that have resisted generations of mathematicians, how can we dismiss the possibility that it will eventually surpass us across all cognitive tasks?  This raises  fundamental questions about human abilities and the prospect of automating them, with enormous implications for research, education, communication---indeed, almost every sphere of human activity.

\subsection{Extinction problem}

 The true origin of that extinction  warning goes back much
earlier, to a 1983 essay by Vernor Vinge, in which he claimed that humans
will soon create ``intelligences greater than our own'' --- that is, a
superintelligent AI (SAI) exceeding human ability across all cognitive
tasks. When this is achieved, Vinge wrote, ``human history will have
reached a kind of singularity, an intellectual transition as impenetrable
as the knotted space-time at the center of a black hole, and the world will
pass far beyond our understanding.''\footnote{Vernor Vinge, \href{https://gwern.net/doc/transhumanism/1983-vinge.pdf}{``First Word,''} \textit{Omni} \textbf{5}, no.~4 (January 1983), p.~10. This is distinct from his later essay, ``The Coming Technological Singularity: How to Survive in the Post-Human Era'' (1993), cited below.}   
The proposed mechanism is recursive self-improvement: humans design an AI capable of improving its own design. The improved system then produces further improved systems at an accelerating rate, leading to a runaway intelligence explosion---the singularity.
Some highly accomplished people invoke precisely this metaphor. On February 23, 2025, Elon Musk declared, ``We are on the event horizon of the singularity.''\footnote{Elon Musk, \href{https://x.com/elonmusk/status/1893810875875889507}{post on X}, 23 February 2025.} In June 2025, Sam Altman went further: ``We are past the event horizon; the takeoff has started.''\footnote{Sam Altman, \href{https://blog.samaltman.com/the-gentle-singularity}{``The Gentle Singularity,''} June 2025.} The event horizon marks a point of no return---a metaphor for the achievement of a system capable of autonomous self-improvement. But has AI actually crossed that threshold and begun the irreversible descent to the singularity?\footnote{For my earlier discussions of these claims, see Sergiu Klainerman, \href{https://www.tabletmag.com/sections/arts-letters/articles/keep-thinking-artificial-intelligence}{``Keep Thinking,''} \textit{Tablet}, 9 January 2026; and \href{https://hxstem.substack.com/p/ai-and-the-doomsday-scenario}{``AI and the doomsday scenario,''} \textit{Heterodox STEM}, 9 August 2026.} 

\subsection{Recursive self-improvement}
``Self-improvement'' covers several different achievements. AlphaZero, for example, learned to play chess and Go through self-play: it improved its performance within a fixed, human-designed learning framework. The Darwin G\"odel Machine goes a step further by repeatedly modifying its coding-agent software and testing the changes on programming benchmarks. Those modifications can also improve its ability to make subsequent modifications. AlphaEvolve proposes and evaluates changes to algorithms, including algorithms used in training its underlying language model.\footnote{Jenny Zhang et al., \href{https://arxiv.org/abs/2505.22954}{``Darwin G\"odel Machine: Open-Ended Evolution of Self-Improving Agents''} (2025); Alexander Novikov et al., \href{https://arxiv.org/abs/2506.13131}{``AlphaEvolve: A coding agent for scientific and algorithmic discovery''} (2025).} These systems can improve not only their performance but also parts of the machinery through which they learn and solve problems. They do not come close to demonstrating a self-sustaining acceleration of general intelligence of the kind envisioned in the singularity scenario.

Biology has been trying for the last seventy years to solve a more modest
problem: the  synthesis  of self-replicating life from ordinary chemistry.  
 The goal is to establish a stable copying mechanism that permits heritable genetic variation, allowing natural selection to act on that variation and generate genuine biological novelty.\footnote{The Miller--Urey experiment, reported in 1953, produced amino acids under conditions intended to mimic the early Earth. The RNA-world hypothesis proposes that RNA, which can both carry genetic information and help chemical reactions occur, preceded today's DNA--protein machinery. The challenge extends beyond making the ingredients: copying must preserve information while allowing variation, and reproduction must work together with a continuing supply of materials and energy. } The results have not just been  disappointing; it appears that the more  
scientists dig into the problem, the more complexity  they    discover, to the point that the synthesis  problem
 appears  today  even  further  away from  resolution.

Just as in the case of AI self-improvement, there has been   substantial progress in producing   specific   components.\footnote{Examples include a \href{https://www.jcvi.org/publications/creation-bacterial-cell-controlled-chemically-synthesized-genome}{chemically synthesized bacterial genome controlling a recipient cell} (2010), \href{https://www.nature.com/articles/s41467-024-53226-0}{Darwinian evolution of self-replicating DNA in synthetic compartments} (2024), \href{https://www.nature.com/articles/s41557-025-01830-y}{RNA replication including a fragment of the RNA catalyst itself} (2025), and \href{https://www.nature.com/articles/s41467-026-69531-9}{integration of DNA replication with membrane-lipid synthesis} (2026). These demonstrate particular functions or combinations of functions; the synthetic-genome experiment uses an existing cell's machinery.} The problem is how to  make their activities sustain one another.
The chemist James Tour poses the challenge in  vivid  terms: grant researchers the requisite molecules, their desired genetic sequences, modern laboratories, and the accumulated scientific literature. The remaining task is how to assemble a functioning cell, and in that respect, he claims, we  are as far as we have ever been from a  solution.\footnote{James Tour, ``An Open Letter to My Colleagues,'' \textit{Inference: International Review of Science}, vol.~3, no.~2 (August 2017), \href{https://inference-review.com/article/an-open-letter-to-my-colleagues}{inference-review.com/article/an-open-letter-to-my-colleagues}. Tour's thought experiment explicitly supplies molecules in the required chemical and stereochemical purities and DNA and RNA in whatever sequences the experimenter desires.}
Biochemist Anthony Futerman's critique also highlights a difficulty directly relevant to laboratory synthesis: reactions that succeed in isolation may fail when other compounds are introduced. In the lipid-synthesis experiment he discusses, adding another compound sharply reduced the yield. Constructing life therefore requires more than producing its components separately: the processes that produce and maintain them must work together without disabling one another.\footnote{Anthony Futerman, ``A Storm in a Primordial Teacup,'' \textit{Inference: International Review of Science}, vol.~7, no.~1 (June 2022), \href{https://doi.org/10.37282/991819.22.37}{doi:10.37282/991819.22.37}.}

Even an optimist like John Sutherland, a leading researcher on the chemical origins of life, acknowledges the difficulty. Writing in 2017 about demonstrating the transition from simple chemicals to a living system, he stated: ``We are not yet close to achieving this end.''\footnote{John D. Sutherland, ``Opinion: Studies on the origin of life---the end of the beginning,'' \textit{Nature Reviews Chemistry} \textbf{1}, article 0012 (18 January 2017), \href{https://doi.org/10.1038/s41570-016-0012}{doi:10.1038/s41570-016-0012}. }

Thus even in biology, where self-replication is an observable fact, constructing a living system from nonliving chemical components remains a formidable challenge. Yet Vinge's prediction is far more ambitious: to build an entirely novel  system, with no known counterpart in the physical world,  that continually enhances its own capacity for self-improvement, eventually surpassing the human mind in every respect.

The analogy concerns autonomy: a living system reproduces and sustains an evolutionary process without scientists continually redesigning it, drawing on materials and energy from its environment. The corresponding ambition for AI is a system that sustains its own improvement without continuing human intervention, drawing only on computing resources, power, and access to data. Improving particular components under human direction\footnote{``Under human direction'' does not mean that each change requires human intervention. Systems such as the Darwin G\"odel Machine can autonomously propose and test changes, including changes that improve their ability to make further modifications. They do so within an experimental framework, with objectives and evaluation procedures established by their designers. Such results do not by themselves demonstrate a self-sustaining acceleration of general intelligence.} is a radically different achievement from establishing such a self-sustaining process.

The prospect of machines taking over invention is central to the singularity scenario. In 1965, I. J. Good argued that the first ultraintelligent machine would be the last invention humans need make.\footnote{I. J. Good, \href{https://doi.org/10.1016/S0065-2458(08)60418-0}{``Speculations Concerning the First Ultraintelligent Machine,''} \textit{Advances in Computers} \textbf{6} (1965), 31--88, section 2. Good qualifies this claim by requiring that the machine remain sufficiently controllable.} Vinge adopted this argument, predicting that superhuman intelligence would create still more intelligent successors on ever shorter timescales.\footnote{Vernor Vinge, \href{https://edoras.sdsu.edu/~vinge/misc/singularity.html}{``The Coming Technological Singularity: How to Survive in the Post-Human Era,''} presented at the VISION-21 Symposium, 30--31 March 1993, section ``What is the Singularity?'' Altman advances a more qualified version of this expectation: AI will accelerate scientific research and help create better AI systems, thereby accelerating further progress. He describes this as an early form of recursive self-improvement, while distinguishing it from fully autonomous self-modification. See Sam Altman, \href{https://blog.samaltman.com/the-gentle-singularity}{``The Gentle Singularity,''} June 2025.} Present AI systems improve performance relative to objectives supplied by us: predicting the next token, minimizing an error, searching for the best fit, writing faster code, or learning through rewards and penalties. But the great discontinuities in science---the changes of paradigm---do not arise merely from moving faster within an existing framework. What is the loss function for discovering algebra, calculus, evolution, genetics, Maxwell's equations, atomic theory, or quantum mechanics? What is the gradient from Newtonian mechanics to general relativity? An optimization process can produce unexpected abilities within a prescribed space of possibilities. It is up to those predicting an intelligence explosion to explain how such a process could generate new conceptual frameworks and sustain their transformation at an ever-accelerating pace.

All this leaves me skeptical of predictions that AGI is around the corner. For the foreseeable future, I expect major improvements in AI requiring new ideas or fundamental changes in design to remain the work of humans rather than AI agents. Yet the machines we build will increasingly outperform us at tasks traditionally regarded as measures of intelligence. This confronts us with a paradox: are the machines we create more intelligent than  us, their creators? Is AlphaZero more intelligent than the engineers who created it simply because it vastly surpasses them at chess or Go? If a powerful AI solves mathematical problems that have resisted generations of mathematicians, does that render mathematicians obsolete?

      \subsection{What is intelligence?}
      
      Before declaring the human mind obsolete, we  better ask what our current  measures of intelligence leave out. Standard IQ tests  do not adequately capture some of our most important capacities: choosing worthwhile goals, sustaining focused effort in pursuing them, and exercising judgment about when to persist, change direction, or abandon them. And, above all,  the human mind has the uncanny ability to refuse to get trapped in any  finite  space; whenever  it becomes aware of its confinement, it  seeks and often   finds a way out of it.
   It can question the assumptions of a given problem, reformulate it, and invent concepts that open possibilities previously beyond its horizon. We are finite creatures endowed with  free will and self-awareness  as well as   an irrepressible intuition of infinity\footnote{As a metaphor for the mind's relation to itself, picture a set containing itself as an element---an image excluded by the axiom of foundation in ZFC, whose repeated unfolding suggests an endless reflection of the mind upon itself.}. Though free will has often been contested,\footnote{See, for example, Derk Pereboom, \href{https://www.cambridge.org/core/books/living-without-free-will/6A4A6C531D58B61C53F6AE4D1A9FD74A}{\textit{Living without Free Will}} (Cambridge University Press, 2001); Sam Harris, \href{https://www.simonandschuster.com/books/Free-Will/Sam-Harris/9781451683400}{\textit{Free Will}} (Free Press, 2012); and Robert M. Sapolsky, \href{https://www.penguinrandomhouse.com/books/592344/determined-by-robert-m-sapolsky/}{\textit{Determined: A Science of Life Without Free Will}} (Penguin Press, 2023).} it is  hard to disagree with  the fact that whenever we   become aware of a
       specific  intellectual  limitation    we   can take  steps to overcome or even  act contrary to the limitation. In that sense   the Socratic injunction ``know thyself'' is both  an act of liberation and engine of  human creativity. 
       
      Humans are slow at performing calculations, but we compensate by thinking in terms of principles and structures that embed the crucial steps of a calculation in the architecture of a theory, making those steps appear natural and easy to verify. Mathematical insight can also precede its expression in words: we may recognize a pattern, picture a geometric relation, or sense a unifying principle before we can formulate it precisely. As Einstein explained to Hadamard, ``The words or the language, as they are written or spoken, do not seem to play any role in my mechanism of thought.''\footnote{Albert Einstein, ``A Testimonial from Professor Einstein,'' Appendix II in Jacques Hadamard, \textit{An Essay on the Psychology of Invention in the Mathematical Field} (Princeton University Press, 1945), pp.~142--143; \href{https://worrydream.com/refs/Hadamard_1945_-_The_psychology_of_invention_in_the_mathematical_field.pdf}{original text}. Einstein goes on to describe thinking through images and seeking words at a later stage.} The written argument gives logical form to ideas whose discovery followed a very different path---one that cannot be reconstructed. This distinction between an idea and its expression is familiar across human creative activities. A composer may hear music inwardly before committing it to notation, just as a mathematician may grasp a structure before finding the words and symbols to describe it. Mozart, for example, distinguished the two stages in a letter to his father about \textit{Idomeneo}: ``the composing is finished, but not the writing out.''\footnote{Wolfgang Amadeus Mozart to Leopold Mozart, Munich, 30 December 1780; \href{https://dme.mozarteum.at/DME/briefe/letter.php?mid=1139}{Digital Mozart Edition, letter and source details}. English translation in \href{https://www.gutenberg.org/cache/epub/5307/pg5307-images.html}{\textit{The Letters of Wolfgang Amadeus Mozart (1769--1791)}}, translated by Lady Wallace.}  

LLMs, by contrast, learn from vast collections of written material, including the mathematical literature available on the internet.\footnote{Some systems are further trained through feedback on their own problem-solving attempts, following the broad principle illustrated by AlphaZero's self-play: improving performance by learning from the outcomes of those attempts.} These texts record ideas after they have found expression in words and symbols, leaving out the thought processes through which a new concept, geometric picture, or unifying scientific principle first emerged---or the work of imagination that gave rise to a  great  work of literature or art. The finished work is available to the machine; the act of imagination that brought it into existence is not.

AI is at its best when human judgment directs it toward clearly defined tasks within restricted domains of knowledge, grounded in verified results or informed expert opinion. Without such direction, it relies on statistical patterns learned from the internet, too often recycling the clichés, muddled opinions, and outright nonsense that abound there. Worse still, it can ``hallucinate,'' inventing facts, quotations, or references and presenting them with all the confidence of an expert.

 Human thought draws on more than language and abstract reasoning: it is shaped by the experience of living in a body. We encounter the world through perception, movement, physical effort, emotion, and interaction with other people.   Our bodies and minds have a complexity, at both microscopic and macroscopic levels, that remains far beyond our understanding.  Giving a machine cameras, sensors, and mechanical limbs does not give it the bodily experience that shapes human thought: hunger, pain, fears and desire---``the heartache and the thousand natural shocks that flesh is heir to.''
  Even the robots envisioned by AI enthusiasts, let alone those available today, fall far short of reproducing the richness of bodily experience which informs our  understanding and creative process.\footnote{Sergiu Klainerman, \href{https://hxstem.substack.com/p/ai-and-the-doomsday-scenario}{``AI and the doomsday scenario,''} \textit{Heterodox STEM}, 9 August 2026.}

Human creativity also draws on taste, wisdom, aesthetic sense, faith, intuition, agency, moral sense, empathy, altruism, ego, competitiveness, and other attributes whose workings do not appear to be reducible to algorithms and calculations. These are inseparable from human intelligence: they shape what we care about, which questions we pursue, and what we judge worth creating.

AI enthusiasts argue that future AI could possess both free will and consciousness, together with  other attributes of the human mind. Some suggest that consciousness could emerge in artificial systems with sufficient complexity and the right organization.\footnote{See Stanislas Dehaene, Hakwan Lau, and Sid Kouider, \href{https://doi.org/10.1126/science.aan8871}{``What is consciousness, and could machines have it?,''} \textit{Science} \textbf{358} (2017), 486--492. Their proposal emphasizes making information globally available within a system and monitoring its own processing. See also David J. Chalmers, \href{https://www.bostonreview.net/articles/could-a-large-language-model-be-conscious/}{``Could a Large Language Model Be Conscious?,''} \textit{Boston Review}, 9 August 2023.} Yet merely invoking ``emergence,'' without specifying the process by which consciousness would arise, is not so different from invoking the hand of God. It simply gives the mystery a fancy, scientific-sounding name rather than a serious explanation.

Whether we regard this prospect as plausible appears to  depend on our view of a deeper question: can the mind, at least in principle, be fully understood in terms of the physical components and physiological processes of the brain? Materialists\footnote{Here is a typical statement  in that vein: `We have already known for a long time from clinical neurology and from my own field, clinical psychiatry, that without a shadow of a doubt there is no aspect of the mind that is not entirely the product of, and utterly dependent on, the physical brain. Disruption, disassembly or enhancement of brain circuitry (subtle or major) can radically alter any aspect of the mind.' Ralph Lewis, \textit{Psychology Today}, \url{https://www.psychologytoday.com/us/blog/finding-purpose/202009/what-actually-is-consciousness-and-how-did-it-evolve}.} answer yes, pointing to the impressive progress of neuroscience in relating subjective experiences to specific patterns of brain activity.  Yet establishing such connections is not the same as explaining the experiences themselves.  Computers provide  a perfect  analogy  for this mind-body problem;  software requires physical hardware to run, but a description of the hardware alone does not explain the program’s logic, structure or meaning.  Likewise a radio can reproduce  the music  produced by a source far away and though the  sound  can be stopped   or more broadly interfered    with     by the circuitry of the  radio, nobody can claim that  the music is  produced by the    radio itself.\footnote{The radio analogy is, of course,  not intended as a theory of the mind. It illustrates a logical distinction: evidence that a capacity depends on a physical mechanism does not, by itself, establish that the mechanism fully explains that capacity. Neuroscience's success in relating mental experiences to brain activity therefore does not alone settle the question of whether the mind is reducible to those activities.}

       I  hold  that objective reality extends beyond what we call the physical world. Mathematics offers the clearest example: its objects are accessible to the mind rather than the senses, yet their properties resist our wishes and constrain our reasoning as firmly as physical objects constrain our actions. As I argued in my essay on Wigner, even our understanding of physical reality depends on mathematical structures whose meaning and coherence extend far beyond their physical applications.\footnote{Sergiu Klainerman, \href{https://inference-review.com/article/reflections-on-an-essay-by-wigner}{``Reflections on an Essay by Wigner,''} \textit{Inference} \textbf{7}, no.~3 (2022), especially the section ``Wigner's Great Mystery and Notions of Reality.''} The mind does more than process information about the material world: it discovers and explores this larger reality. 

\subsection{AI and mathematics}

Even though AI may not acquire consciousness or surpass us at every task anytime soon, it is already having a major impact on society. As a microcosm of that impact, consider mathematics, my own field of research and teaching. Paradoxically, although AI is itself largely a product of mathematics\footnote{Mathematics provides the fundamental conceptual framework of modern AI. Neural networks, their learning rules, and the objectives that guide their training are mathematical constructions, drawing on linear algebra, calculus, probability, statistics, and optimization. Computing power, data, and engineering turn these mathematical constructions into AI systems with billions of parameters, trained on enormous amounts of data.}, there is no area of science in which its impact has been greater. In the past month alone, AI-generated work has claimed breakthroughs on a large number of famous mathematical problems, including work on four of the Clay Millennium Prize Problems: the Navier--Stokes existence and smoothness problem, the Riemann hypothesis, the Hodge conjecture, and the Birch and Swinnerton-Dyer conjecture. These claimed results, part of the October 6 mass release described above, have shaken the mathematical community to its core.  Even before these announcements, Jacob Tsimerman, one of the 2026 Fields Medalists, concluded his lecture at the International Congress of Mathematicians in Philadelphia that July by predicting that AI would become superhuman at mathematics within a few years, producing proofs and proposing new conjectures and insightful definitions.\footnote{Marianne Freiberger and Rachel Thomas, \href{https://plus.maths.org/icm-2026-postcard-philadelphia}{``ICM 2026: A postcard from Philadelphia,''} \textit{Plus}, 5 August 2026. Their eyewitness account reports Tsimerman's remarks at the conclusion of his Fields Medal lecture.} These predictions are having a profound effect on mathematicians, especially younger researchers, many of whom wonder whether it still makes sense to pursue a career in mathematics.

 Though the challenges posed by AI's extraordinary recent advances are serious, Tsimerman's predictions, delivered with full confidence and the authority of a newly awarded Fields Medal, are, in my view, as irresponsible as predictions that AGI will soon kill us all---particularly given their demoralizing effect on young mathematicians. We need to distinguish serious concerns about how AI is presently used from sweeping predictions presented with unwarranted certainty.

The significance of these announced AI breakthroughs depends first on whether the claimed proofs dumped by AI are correct.\footnote{Formal verification in a proof assistant such as Lean is a powerful aid to mathematical verification, but the label ``Lean-verified'' should not by itself be treated as conclusive validation of an announced result. One must examine the precise statement formalized, its assumptions, and whether the proof is complete, and establish that these faithfully capture the intended mathematical claim. Lean can substantially strengthen this human assessment by checking the formal deductions. It does not replace the work of understanding what has been proved, how it relates to earlier mathematics, or why it matters. See \href{https://lean-lang.org/theorem_proving_in_lean4/Introduction/}{\textit{Theorem Proving in Lean 4}}, Introduction. See also Alexander Bastounis, Fabian Circelli, and Anders C. Hansen, \href{https://arxiv.org/abs/2610.08144}{``Navier--Stokes lost in translation: Why Lean verification of AI autoformalisation does not guarantee correct natural language proofs,''} arXiv:2610.08144, 6 October 2026. The authors identify discrepancies between the announced Navier--Stokes argument and its Lean formalisation, arguing that verification of the latter does not establish the correctness of the former.} But even if the claimed results are correct---and many are---they do not by themselves justify Tsimerman's predictions. A proof may resolve a major conjecture through a new conceptual idea, a powerful extension of existing methods, or computations previously beyond reach. These are all achievements, but they provide different grounds for assessing predictions about the future of mathematics.

This dim view of the future of mathematics, in which humans are replaced by machines, rests on the mistaken belief that the main purpose of mathematics is to prove theorems. But a mathematician does more than turn coffee into theorems, as R\'enyi's famous quip, often repeated by Erd\H{o}s, would have it. Like any other science, mathematics seeks to discover significant truths and understand their connections. In the natural sciences, claims about the physical world are tested through observation and experiment and unified by mathematical theories that yield further testable predictions. In mathematics, proof is the means by which discoveries are validated. If AI can help us find proofs of statements that have resisted our best efforts, so much the better. The natural sciences already rely extensively on computer-assisted calculations. Automating parts of mathematical proof would bring mathematics closer to their situation: machines perform more of the technical work, while scientists seek the ideas that give that work meaning and direction.

Mathematicians have long considered the possibility of automating proofs. Hilbert's tenth problem, posed in 1900, asked for an algorithm to determine whether a polynomial equation with integer coefficients has an integer solution. His later \textit{Entscheidungsproblem}, or decision problem, asked for a procedure that would always terminate and determine whether any given statement in first-order logic\footnote{First-order logic is a formal language in which quantifiers such as ``for every'' and ``there exists'' range over individual objects in a chosen domain. In arithmetic, these objects are numbers. Second-order logic also allows quantification over properties or relations on those objects, as in ``for every property of natural numbers,'' which permits induction to be expressed as a single axiom. First-order arithmetic instead expresses induction through an axiom schema, with an instance for each formula. First-order set theory can nevertheless describe sets, functions, and relations by treating them as objects in its domain.} is provable.

Church and Turing showed in 1936 that no such decision procedure exists. This does not rule out automated proof search: within a fixed, effectively specified formal system, an algorithm can systematically examine all possible proofs and eventually find a proof whenever one exists. If none exists, however, the search may continue forever. Searching through all possible proofs is usually hopelessly slow. AI can make the search more effective by suggesting promising approaches and useful intermediate steps, then changing course when an attempt fails. This can bring previously inaccessible proofs within reach, but it does not overcome the logical limitations established by Church and Turing: there is still no procedure guaranteed to settle every problem. More importantly, its growing ability to find proofs must also be judged against the deeper goal of mathematical understanding. If a proof succeeds only because a machine can examine a myriad of cases faster than humans, without revealing any unifying principle, it may establish a fact while offering little insight into why it holds. Good mathematical proofs do more than validate statements: they reveal connections and, at their best, introduce conceptual frameworks that lead to further discoveries.

A striking example is the contrast between Abel's proof that the general polynomial equation of degree five or higher cannot be solved by radicals and Galois's subsequent theory. Abel settled the impossibility question; Galois revealed the underlying structure, relating solvability to the group of symmetries of the roots. Abel's original proof has been essentially forgotten, while Galois's conceptual framework transformed algebra and became fundamental throughout mathematics and science.

The uniformization theorem provides another striking  example. Its PDE proof, based on finding a conformally related metric of constant curvature, is more revealing than the classical complex-analytic proof, going back to Riemann. It exposes the relation between conformal geometry, curvature, and nonlinear differential equations, and points toward higher-dimensional generalizations, notably the Yamabe problem of finding metrics of constant scalar curvature within a given conformal class.
The celebrated Hamilton--Perelman proof of the Poincar\'e conjecture carries this geometric PDE approach much further. Hamilton introduced Ricci flow, an evolution of the metric governed by a nonlinear PDE, with the aim of obtaining canonical geometries; Perelman's  later breakthroughs made it possible to establish the  broader geometrization conjecture. 

In my own work on the stability of black holes, we now have a complete proof spread across many highly technical papers totaling well over two thousand pages.\footnote{There are two approaches to nonlinear stability of Kerr black holes in the full subextremal range. One combines my earlier work with J\'er\'emie Szeftel, our joint work with Elena Giorgi, and work by Dawei Shen with recent results of Szeftel and of Siyuan Ma and Szeftel; see Szeftel, \href{https://arxiv.org/abs/2609.40068}{``Kerr stability in the full subextremal range,''} and the references therein. Peter Hintz announced an independent result using microlocal methods developed in his earlier work with Andr\'as Vasy; see \href{https://arxiv.org/abs/2606.28253}{``Nonlinear stability of subextremal Kerr black holes.''} His  result requires significantly stronger assumptions on the asymptotic behavior of the initial data.} I can think of no more important goal than finding a shorter, more conceptual proof---one that reveals the essential mechanisms and opens the way to further developments.

Mathematics abounds in examples showing that the ideas revealed by a proof are more important than the result it establishes.
 Nothing illustrates this more dramatically than new methods that transform mathematics even when, at first glance, they do nothing more than prove old theorems. Cartesian geometry provided a systematic way to translate problems in Euclidean geometry into algebraic equations. Its significance extended far beyond new proofs of familiar results: it transformed the relation between geometry and algebra, helped pave the way for calculus, and provided a mathematical language for the emerging science of physics.  Grothendieck, one of the twentieth century's great builders of mathematical theories, compared two approaches to proving a theorem with two ways of opening a nut: striking it with a hammer and chisel, or immersing it in water until the shell softens enough to open by hand. The second captures his approach to mathematics: patiently developing concepts and structures until a previously resistant problem  becomes a simple natural consequence  of a broader  theory.\footnote{Alexandre Grothendieck, \textit{R\'ecoltes et Semailles}, section 122, ``La mer qui monte\ldots''; see the passage reproduced by the \href{https://csg.igrothendieck.org/quotes/}{Grothendieck Institute}.} He also described his approach through the image of a rising sea: ``Little by little it becomes a peninsula, then an island, then an islet, which in turn is finally submerged.''\footnote{Grothendieck, \textit{R\'ecoltes et Semailles}, section 122, ``La mer qui monte\ldots''; translated from the \href{https://geotoppos.org/divers/ReS.pdf\#page=686}{French original}. Here ``it'' refers to the resistant ground gradually surrounded by the sea.}
 
Physicists and engineers use mathematical theories---their concepts, equations, and methods. They rarely concern themselves with the precise statements of theorems, much less their proofs. The extraordinary usefulness of mathematics rests primarily on its power to provide  theoretical frameworks for understanding and describing the world. When Einstein wrote in ``On the Method of Theoretical Physics'' that ``the truly creative principle resides in mathematics,'' he was pointing to this power of mathematics to guide physical discovery, rather than to any particular mathematical result.\footnote{Albert Einstein, \href{https://doi.org/10.1086/286316}{``On the Method of Theoretical Physics,''} \textit{Philosophy of Science} \textbf{1}, no.~2 (April 1934), 163--169, p.~167. The Herbert Spencer Lecture, delivered at Oxford on 10 June 1933.} He was keenly aware that the search for his beloved unified theory might require mathematics that had yet to be discovered.

An obsessive fixation on solving famous problems can sometimes be counterproductive. The announced AI solution to the forced Navier--Stokes Clay problem\footnote{The announced construction concerns the far more questionable case of the Navier--Stokes equations with a specially designed external forcing term, rather than the unforced initial-value problem. Forced counterexamples are explicitly permitted by the Clay problem's formulation. For my assessment of their physical relevance, see Sergiu Klainerman, \href{https://inference-review.com/article/ai-will-not-make-mathematicians-obsolete}{``AI will not make mathematicians obsolete,''} \textit{Inference} \textbf{9}, no.~2 (September 2026); also published in \href{https://hxstem.substack.com/p/ai-will-not-make-mathematicians-obsolete}{\textit{Heterodox STEM}}, 14 September 2026.}, for example, while technically impressive, teaches us little about the behavior of real fluids. The exceedingly complicated blow-up solutions it constructs are highly unphysical and most probably highly unstable, though this remains to be established. Before computational assistance became available, I would have regarded encouraging someone to construct such a solution by hand as professional malpractice. By contrast, showing that solutions remain regular for all time outside a suitably defined exceptional set of initial conditions is a foundational question---a first step in a program aimed at understanding the behavior of fluids, including turbulence. That problem remains far beyond the current abilities of both humans and AI.

   Another misconception is that mathematics consists of fixed fields---such as number theory, algebraic geometry, topology, differential geometry, and probability---organized around lists of open problems. Tsimerman refers to mathematics as a closed system. But it appears closed only if we reduce its purpose to settling  famous problems within these established fields. Such a setting is particularly favorable to AI: the questions are already formulated, the relevant concepts developed, and a substantial body of methods and results available to draw upon. AI can search this literature, combine techniques, and test possible routes to a proof. What this picture leaves out is the creation of new concepts, new questions, and new fields of mathematics.

In fact mathematics continually enlarges and transforms its frameworks. Contact with other sciences, especially physics, has been indispensable to this renewal. Such exchanges were central to the golden age stretching from  Descartes,  Galileo and Newton to  Gauss, Riemann, Maxwell and Einstein. As Poincar\'e emphasized, nature presents questions that mathematical imagination alone would never have conceived, preventing mathematics from endlessly turning in circles.\footnote{Henri Poincar\'e, \textit{The Value of Science}, chapter V, ``Analysis and Physics,'' section III; English text published in \textit{Popular Science Monthly}, vol.~70 (February 1907), \href{https://en.wikisource.org/wiki/Popular_Science_Monthly/Volume_70/February_1907/The_Value_of_Science:_Analysis_and_Physics_VI}{full text}. The passage reads: ``The combinations which numbers and symbols may form are an infinite multitude. In this multitude how shall we choose those which are worthy to fix our attention? Shall we let ourselves be guided solely by our caprice? This caprice, which itself would besides soon tire, would doubtless carry us very far apart and we should quickly cease to understand each other. But this is only the smaller side of the question. Physics will doubtless prevent our straying, but it will also preserve us from a danger much more formidable; it will prevent our ceaselessly going around in the same circle. History proves that physics has not only forced us to choose among problems which came in a crowd; it has imposed upon us such as we should without it never have dreamed of. However varied may be the imagination of man, nature is still a thousand times richer. To follow her we must take ways we have neglected, and these paths lead us often to summits whence we discover new countries. What could be more useful!''} AI's growing ability to resolve conjectures within established fields may encourage a renewal of this interaction. By freeing mathematicians to concentrate on conceptual development, it may stimulate both new unifying theories within existing fields and the creation of new fields in response to questions arising in other sciences.

Even if theorem proving were eventually delegated entirely to AI---a prospect still far from today's reality---mathematics would retain its central creative tasks: deciding which questions matter and developing the concepts and theories that make their answers illuminating. This would mark not the end of mathematics, as many fear, but the beginning of an exciting new stage in its development.

\subsection{There's never been a better time to go into math}

Grant Sanderson, creator of 3Blue1Brown, chose this provocative title for a recent video arguing that AI's advances open new opportunities for mathematicians.\footnote{Grant Sanderson, \href{https://www.youtube.com/watch?v=LKrwhjOI8Nc}{``There's never been a better time to go into math,''} 3Blue1Brown, 16 October 2026.} I share his optimism: mathematics extends far beyond proving theorems, and automating parts of that work may indeed create room for new ideas and directions. But realizing these opportunities will require confronting serious challenges---in how we assign credit, train young mathematicians, and preserve the open exchange of ideas.\footnote{For a wide range of mathematicians' reactions to OpenAI's release, see \href{https://proofsandprompts.com/2026/10/08/100-reactions-to-100-solutions/comment-page-1/}{``100+ reactions to 100+ solutions,''} \textit{Proofs and Prompts}, 8 October 2026. Contributions range from Ryan Alweiss's enthusiasm for the opportunities opened by collaboration with AI to Tristan Humbert's concerns about the disruption of young researchers' work and the standards of the released manuscripts.}

My own work with AI illustrates a new mode of mathematical research: a constant dialogue guided by my goals, ideas, and judgment. AI helps explain papers, check calculations, develop approaches I suggest, and test examples and counterexamples. It also assists with LaTeX and exposition. Through repeated questioning and correction, I examine its arguments and assumptions: a plausible answer may contain an error or leave the central difficulty untouched. Making each step explicit helps clarify the obstruction and decide how to proceed.

If OpenAI releases the internal model behind its recent collection of mathematical results, I may also be able to ask AI to carry an approach through to a complete proof.   If the  proof    strikes me as interesting,    it  could open a further phase of research: digesting the proof, identifying its essential ideas, and developing them into concepts and methods that reach beyond the original problem. The proof would then become a starting point for new questions and a broader theory. I would welcome this: it would allow me to concentrate on the part of mathematics I enjoy most---finding the concepts that unify what I know and developing a vision for new directions of research.

All this is very promising, but there are also serious risks. AI could greatly expand what mathematicians can accomplish while undermining the community that makes those accomplishments possible. We must distinguish between challenges that are serious but manageable and those that threaten the continued vitality of mathematics as a human enterprise.

OpenAI's hostile mass release exposes a darker side of AI. Hundreds of proposed results arrive at once, with little explanation of the central new ideas that make the arguments plausible or of how they build on the work of mathematicians in the field. The burden of checking the proofs, identifying those ideas, and understanding their implications falls heavily on the mathematical community. AI can generate papers much faster than mathematicians can check and understand them. The resulting imbalance threatens to overwhelm the attention, rigor, and conceptual understanding on which mathematical progress depends.

Hugo Duminil-Copin, whose reaction was cited earlier, described OpenAI's release as having ``annihilated'' his field of research. In his October 14 Numberphile interview, ``The Strip Mining of Mathematical Proofs,'' he explains what such a disruption means for the mathematicians who have devoted their lives to developing the field.\footnote{Hugo Duminil-Copin, interview with Brady Haran, \href{https://www.numberphile.com/podcast/hugo-ai}{``The Strip Mining of Mathematical Proofs,''} \textit{The Numberphile Podcast}, 14 October 2026.} He has also warned that young researchers will be especially affected and that fear of companies using publicly shared ideas to claim results threatens the culture of open exchange.\footnote{David Larousserie, \href{https://www.lemonde.fr/en/science/article/2026/10/12/openai-s-latest-blow-to-mathematicians-they-have-annihilated-my-field-of-research_6758515_10.html}{``OpenAI's latest blow to mathematicians: `They have annihilated my field of research',''} \textit{Le Monde}, 12 October 2026, reporting his remarks at the Soci\'et\'e Math\'ematique de France on 8 October and a subsequent interview.} The mining metaphor captures the danger of treating mathematical ideas as deposits from which finished results can be extracted, without sustaining the community that develops those ideas. That community does more than supply problems: it cultivates understanding, trains successors, and creates new directions. An abundance of extracted results would be a poor exchange for damage to those activities.

Tristan Buckmaster's account is particularly telling. He explicitly credits Diego C\'ordoba and Luis Mart\'inez-Zoroa with the ideas behind the program that he and Levent Alp\"oge developed with AI assistance. He also describes pressure to release their results before they had time to turn the generated proofs into readable mathematical arguments. In recounting his exchanges with OpenAI's S\'ebastien Bubeck, Buckmaster also describes an implicit threat to his career: when he said he would make the circumstances public if OpenAI released its result in the proposed manner, he reports being asked, ``Why would you ruin your career?'' When he challenged that remark, the response, according to his statement, was: ``If you don't want me to be nice, then I don't have to be nice.''\footnote{Tristan Buckmaster, \href{https://cims.nyu.edu/~tristanb/statement.pdf}{public statement on his work with Levent Alp\"oge and their interactions with OpenAI}, September 2026.}

Whether OpenAI intended these releases to generate publicity, justify its enormous valuations, or bolster its claim to be approaching superintelligence, we do not know. But the announcements have been deeply demoralizing, especially for young researchers, some of whom now question whether they have a future in mathematics.

Flooding mathematics with claims faster than they can be assessed, pressuring researchers to release work before it is ready, and threatening those who object are abuses of power that the mathematical community must confront together. The mathematical community, and more broadly the entire scientific community, must establish enforceable standards of verification, attribution, and open exchange, and make compliance a condition of collaboration with AI companies.

The only acceptable model is one in which AI works under the direction and supervision of the researchers responsible for the work. They set the goals, judge the significance of the results, and take responsibility for the work.

But even under this model, serious difficulties remain. First, if providers such as OpenAI or Anthropic can use a researcher's exchanges with AI to train their models, unpublished ideas may find their way into responses given to others. Researchers need firm guarantees that their work will remain confidential and will not be used for training without their explicit consent.

Second, even if confidentiality is secured, the speed with which AI can extend and adapt results threatens to disrupt both mathematical communication and the allocation of academic credit. An idea shared in a lecture or preliminary manuscript could be developed and published by others before its originator has time to carry it through. If credit goes chiefly to whoever first produces a complete proof, mathematicians may become reluctant to share the ideas on which collective progress depends.

Third, collaboration with AI is producing a rapidly growing volume of papers that threatens to overwhelm journals and the peer-review system.\footnote{For an example of the pressures associated with high-volume AI-assisted research, see Gabriel Sison, \href{https://www.thecollegefix.com/publisher-removes-200-papers-university-of-chicago-statistician-wrote-in-under-1-year/}{``Publisher removes 200 papers University of Chicago statistician wrote in under 1 year,''} \textit{The College Fix}, 12 October 2026. Citing reporting by \textit{The Washington Post}, the article states that SSRN removed 253 papers written or co-written by Nicholas Polson within the preceding year, quoting SSRN's explanation: ``The most significant issue was the unusual submission activity.'' SSRN is a preprint platform; these removals should not be confused with retractions of peer-reviewed journal articles. Polson acknowledged using AI and defended the quality of his work.} Generating manuscripts can take far less time than checking their correctness, identifying their original contributions, and assessing their significance. The resulting flood risks burying important ideas beneath work that few mathematicians have time to read carefully. 

Fourth, AI challenges mathematical training based on written homework, though greater reliance on oral examinations can help address this. The more serious difficulty concerns doctoral students, whose initial research training typically consists of understanding and extending work done by others. Students can now ask AI to carry out much of this work, potentially obtaining results without developing the understanding and independence that the  doctoral training  is meant to cultivate.

Protecting unpublished work entrusted to AI requires resources and institutional guarantees beyond the reach of the mathematical community alone.

I therefore propose an academic AI consortium bringing together the NSF, American universities, the Simons Foundation, and other private institutions. Its mission would be to provide powerful AI tools dedicated to scientific research, with enforceable guarantees that unpublished ideas will neither be used for training nor shared with others without the researchers' explicit consent. The consortium would negotiate protected access to commercial models, operate models itself, or support their development. Scientists should not have to choose between access to the best research tools and control over their own ideas. Creating this consortium should be a top priority for the entire scientific community. 

The remaining challenges are ours, as mathematicians, to meet: preserving the open exchange of ideas, developing better measures of success that recognize broad intellectual contributions, establishing high ethical standards and effective means of enforcing them, training mathematicians capable of independent thought, and finding new ways to teach our beloved subject. None of these tasks will be easy, but I see no reason to doubt that we can meet them. That is why I share Sanderson's optimism about the opportunities awaiting the next generation of mathematicians. 

\subsection{A few  comments on the broader impact of AI on society}

Public discussion of AI spans a broad spectrum, from Gregory Conti's deep pessimism to Tyler Cowen's optimism tempered by concern and Shyam Sankar's forceful rejection of AI alarmism.\footnote{Gregory Conti, \href{https://www.compactmag.com/article/big-techs-war-on-human-achievement/}{``Big Tech's War on Human Achievement,''} \textit{Compact}, 27 July 2026; Tyler Cowen and Avital Balwit, \href{https://www.thefp.com/p/ai-will-change-what-it-is-to-be-human}{``AI Will Change What It Is to Be Human. Are We Ready?,''} \textit{The Free Press}, 12 May 2025; Shyam Sankar, \href{https://www.thefp.com/p/dangerous-ideology-effective-altruism-artificial-intelligence}{``The Dangerous Ideology Behind the AI Warnings,''} \textit{The Free Press}, 14 September 2026.} These disagreements extend beyond the widely exaggerated fears of human extinction to what AI's growing capabilities will mean for human creativity, work, and the possibility of a meaningful life.

Cowen and Balwit welcome AI's promise while worrying about our place in a world where machines surpass us. But that worry already concedes the premise I have challenged throughout this essay.

 My colleague  Gregory Conti fears that AI will destroy the opportunities for meaningful achievement through which people develop their abilities and find a vocation. His concern goes beyond the loss of employment: it concerns the loss of reasons to strive. But his conclusion assumes that machines will exhaust the possibilities for worthwhile human endeavor. The discussion of mathematics above suggests otherwise. A problem solved can reveal connections we had not suspected and bring into view questions we could not previously formulate. AI may further blunt the intellects of the lazy, as Conti rightfully fears, but it may also inspire the ambitious and the daring to reach new heights of scientific achievement or artistic expression and deepen our understanding of the human condition.

It is not surprising that  academics, in particular,  feel threatened. Earlier technological revolutions were mostly experienced by educated professionals as disruptions to other people's jobs, chiefly those requiring less formal education. AI brings that disruption directly into their own world. Activities long regarded as protected by advanced education---writing, programming, research, and professional judgment---are increasingly open to automation. Some of the alarm surrounding AI may reflect this reversal: those accustomed to explaining why others must adapt now face that demand themselves.

In many ways, I think that the opposite of what Conti fears may come to pass. If capable robots become plentiful and inexpensive, their combination with AI could give ordinary people access to assistance once reserved for aristocrats. Much of the labor required to maintain everyday life could be delegated to machines. Why assume that the resulting freedom must impoverish our lives? Less compulsory work could mean more time for family, friendship, learning, artistic creation, and public service. Leisure does not guarantee a meaningful life, but neither does economic necessity. The opportunity would be extraordinary, provided that access to these benefits is widely shared.

The dangers that deserve our immediate attention concern how people design, deploy, and use these systems: to concentrate power, manipulate others, evade responsibility, or exercise force. Errors and failures of control also matter, even without malicious intent. None of these dangers requires a machine with ambitions of its own. They require us to take responsibility for the powers we are creating.

Sankar identifies effective altruism as the ideology behind much of the alarm.\footnote{Effective altruism seeks to identify how to do the most good through evidence and reasoning, under the presumption that long-term consequences can somehow be calculated. Communist eschatology rested on a similar pretension: that humanity's future could be scientifically foretold.} His central objection is that speculative claims about humanity's distant future are being used to override the interests and freedoms of people living today. As he puts it, ``Putting disputed moral judgments and uncertain prophecies into an equation does not turn them into scientific facts.'' He compares those claiming authority on this basis to Lenin's vanguard: a small group convinced that its privileged understanding of humanity's future entitles it to direct everyone else. Against this ambition, he points to concrete uses of AI in hospitals and shipbuilding that expand what people can accomplish.

One of the most frequently voiced fears concerns the so-called alignment problem: how to ensure that AI acts in accordance with human intentions and values. But human beings themselves disagree profoundly about values, purposes, and acceptable conduct. With whose values should AI be aligned, and who gets to decide? We should demand reliable behavior, clear limits, and accountability in particular applications. Expecting a machine to resolve humanity's moral disagreements is an absurd  ambition.

Ensuring that AI performs its assigned tasks reliably is a practical engineering challenge. The fear that a superintelligent system might escape human control is a different concern, arising from the ambition to build such a system in the first place.

Paradoxically, the main AI companies warn that we may lose control of AI while racing to build machines intended to surpass us in every intellectual capacity. We should instead build powerful AI tools to advance science, transform medicine, and strengthen industry. Our task is to enlarge human possibilities, not to chase the fantasy of making ourselves obsolete.

\section*{Notes}
\printendnotes
\end{document}